\documentclass[12pt]{article}

\newcommand{\invisible}[1]{}
\newcommand{\numero}[1]{\begin{center}
{\bf \large #1 } \end{center}}

\begin{document}

\section*{Set-theoretical mathematics in Coq} 

\noindent
Carlos Simpson\newline
carlos@math.unice.fr\newline
CNRS, Laboratoire J.A. Dieudonne\newline
Universite de Nice-Sophia Antipolis

{\small {\bf Abstract:}
We give a brief discussion of some of the issues which have arisen in the 
course of formalizing some classical set-theoretical mathematics in the Coq system.
This sprouts from, expands and replaces a chapter of math.HO/0311260 which will be removed in
revision, and also contains as a tar-attachment to the source file the revised and
expanded version of the proof development which had been attached to math.HO/0311260. 
}

\numero{Introduction}
This is an expanded version of one of the chapters of \cite{ctpm} which
gave some details about the development which was attached to that source file.
The present preprint also comes with a development attached to its 
source file,
a modified version of the one of \cite{ctpm}: for one thing, it has
been translated to Coq Version 8 syntax (with the help of the automated translation
tools attached to Coq v8); the files have been reorganized and renamed; and we have
plunged further into the notions of well-ordered set, getting to the theorem that any
set can be well-ordered, ordinals, the isomorphism between a well-ordered set and its
corresponding ordinal, the definition of cardinals and the Bernstein-Cantor-Schroeder theorem
(which we deduce as an easy consequence of the theory of ordinals). 

All of the mathematics which we cover here has already been formalized many times over in
different proof systems, such as Mizar, Isabelle/ZF, and more recently Metamath, etc...
As far as I know, it hasn't been fully done in Coq although many parts have been. 
But that is perhaps an uninteresting distinction.  
The present note should be thought of as more of an experiment in axioms, semantics, style, 
notation and organization. 
As justification
for why to write it up, apart from the fact that you have to start somewhere, 
the real point is that it is important to try a variety of different
points of view on the same subject,
and it is important that each of us let everybody else know about these attempts. 

In order fully to grasp the wide range of problems which are encountered
in formalizing mathematics, we need to see what it looks like ``from the inside'', that is to say
from the perspective of document creation. In the textual part of this preprint, I will try to
make remarks in this direction.

Because of the only moderate interest of adding another formalization to those which are already
available, we will in the accompanying note (i.e. the text below) try to describe simply and briefly
the basic outline of what is being done without entering into excessive detail. 
We refer to the Coq reference manual for all explanation of how the system works, and
eventually to our files themselves (or, more usefully, to any of the very numerous examples
the reader will find on the web) for examples. Also we refer to \cite{ctpm} for a more
complete list of references. 

For information, here are the files in the tar archive connected with this source file and 
which make up the proof development:
\begin{verbatim}
axioms.v  tactics.v  set_theory.v  functions.v  notation.v
order.v transfinite.v ordinal.v cardinal.v algebra.v topology.v
\end{verbatim}
We have also included the files {\tt qua.v} and {\tt h1.v} which are refered to in \cite{ctpm}.
These do not constitute a part of the rest of the proof development and compile separately.
They are translated to Version 8.0 syntax using the Coq translation tools,
but otherwise are the same as were attached to 
\cite{ctpm}. 

\numero{Axioms}

We start by
printing the parameters and axioms at the start of our file \verb}axioms.v}:
{\scriptsize
\begin{verbatim}
(*** interpret types as being classical sets ***)
Definition E := Type.

(*** elements of a set are themselves sets ***)
Parameter R : forall x : E, x -> E.  
Axiom R_inj : forall (x : E) (a b : x), R a = R b -> a = b.

Definition inc (x y : E) := exists a : y, R a = x.
Definition sub (a b : E) := forall x : E, inc x a -> inc x b.

(*** a set is determined by its elements ****)
Axiom extensionality : forall a b : E, sub a b -> sub b a -> a = b.

(*** extensionality for general product types ***)
Axiom prod_extensionality : forall (x : Type) (y : x -> Type) (u v : forall a : x, y a),
(forall a : x, u a = v a) -> u = v.

(*** the axiom of choice ****)
Inductive nonemptyT (t : Type) : Prop := nonemptyT_intro : t -> nonemptyT t.
Parameter chooseT : forall (t : Type) (p : t -> Prop), nonemptyT t -> t.
Axiom chooseT_pr : forall (t : Type) (p : t -> Prop) (ne : nonemptyT t), ex p -> p (chooseT p ne).

(*** the replacement axiom: images of a set are again sets *********)
Parameter IM : forall x : E, (x -> E) -> E.
Axiom IM_exists : forall (x : E) (f : x -> E) (y : E), inc y (IM f) -> exists a : x, f a = y.
Axiom IM_inc : forall (x : E) (f : x -> E) (y : E), (exists a : x, f a = y) -> inc y (IM f). 

(*** the following follow from the above but are written as axioms for clarity ********)
Axiom excluded_middle : forall P : Prop, ~ ~ P -> P.
Axiom proof_irrelevance : forall (P : Prop) (q p : P), p = q.

(*** equivalent propositions are equal **********************)
Axiom iff_eq : forall P Q : Prop, (P -> Q) -> (Q -> P) -> P = Q.

(*** the elements of nat are realized as the standard finite ordinals ***)
Axiom nat_realization_O : forall x : E, ~ inc x (R 0). 
Axiom nat_realization_S : forall (n : nat) (x : E), inc x (R (S n)) = (inc x (R n) \/ x = R n). 
\end{verbatim}
}

The principal remark about these axioms is that we are combining type-theory with
set-theory in that the type \verb}E} of sets is declared to be same as \verb}Type}.
This was called ``Option 10.4'' in \cite{ctpm}. It is intended to 
open up the possibility of using the advantages of type theory and set theory alternatively,
whichever is the most convenient for each part of the argument. As a matter of personal taste,
we stick to set-theoretic conventions whenever possible. Given that, one of the main places
where it becomes interesting to change into the type-theory world is to make use of Coq's 
inductive data structures. These are used for natural numbers (and later, integers will also be used),
since in that case we have access to an excellent standard library on Peano arithmetic.
The next main place where this is used is in setting up a system of notation; this will be discussed
more extensively below. 

The \verb}iff_eq} axiom turns out to be extremely useful in structuring the document.
This observation comes from C. Raffali's expos\'e at Nice, see \cite{Phox}. A lemma whose conclusion is
an equivalence between two propositions, is stated as giving the equality of these propositions.
Rewriting can then be used to replace one statement by the other, either in the goal or in a hypothesis
of the context. This very conveniently replaces the \verb}<->} notation which is not sufficiently
supported in Coq. This phenomenon could be system-dependent: with sufficient support (i.e. rewriting-like
tactics) it might well be enough to make due with \verb}<->}. Under the philosophical circumstances of
the present work there is no need to have qualms about using \verb}iff_eq}.

Another important operation which is not sufficiently supported in Coq is the logical ``and''.
We introduce a notation \verb}&} at the highest priority level, permitting its use without
too many parentheses. We also introduce some tactics designed to decompose conjunctions which occur either
in the hypotheses of the context, or in the goal. 

\numero{Tactic notation}

We introduce shorthand notation for most of the commonly used tactics as well as for our own tactics.
Whenever possible, tactics consist of two lowercase characters. For example the 
conjunction-decomposition tactic mentionned above is \verb}ee}. The assumption tactic becomes
\verb}am}, the trivial tactic \verb}tv}, the apply tactic \verb}ap}, and
rewrite tactics \verb}rw} or \verb}wr} depending on the direction,
with \verb}rwi} and \verb}wri} for rewriting in a hypothesis. Proofs start with the intros tactic \verb}ir},
and the tactic for unfolding the head constructor of a hypothesis is \verb}uh}.  

This babble is designed to facilitate typing
proof scripts. Curiously enough, it doesn't diminish readability of the proofs, and might even 
contribute positively in this direction. This is because the tactics which are shorthanded are the 
standard type of proof steps which are skipped entirely in any natural-language proof. Thus what
stands out in a proof are the longer expressions, which are exactly those in which the user has to 
supply something by hand. This includes the names of previous lemmas which are used, and the expressions
of terms which are introduced in the course of the proof. This is the main information which the
reader needs in order to reconstruct in his mind how the proof is being attacked. 

Again this idea came out of Raffali's visit. Indeed, in Phox, the tactics which are abbreviated here to two 
characters, disappear almost entirely. It might be nice to see if this could be done in Coq; however
it seems that Raffali fine-tunes things using his notion of ever-modifiable tactics, so that in a certain
sense the tactics which will enter into a proof are specified beforehand. In our present context of 
trying to type proofs as fast as possible, it seems to be a reasonable compromise to have what is
certainly a wide array of tactics, but which are stationary and can be typed very quickly. 

\numero{Choice and function definition}

One of the characteristics of set-theoretic mathematics which we would like to maintain is
the non-constructive nature, not that this is a goal in and of itself, but because it allows us
to take a very untyped approach to definitions of things. This goes hand in hand, in a philosophical
sense at least, with the axiom of choice: if you know that something exists, then choose one without
worrying about how to construct it explicitly or even how to specify it uniquely. To go even further:
if you are not sure that it exists, choose something which is what you want if one exists, and is
something else if there are none. This is Hilbert's $\epsilon$-function. It has a corollary for
the writing style, which is that function evaluation is defined everywhere. In fact, all of the constructors
which appear in the axioms are defined on all sets (i.e. they don't require a predicate to hold in order
to give a value), so automatically anything else you can define also  has this property. In particular,
evaluation of a function \verb}f} at an element \verb}x} is done by choosing a pair $(x,y)\in f$ (if it exists)
and evaluating to the second projection \verb}y} of this chosen pair. Thus the evaluation 
which we denote by \verb}V x f} is defined for any two sets \verb}f} and \verb}x}.  

The combination of everywhere-definedness and untyped-style
was suggested by L. Lamport in \cite{Lamport} \cite{LamportPaulson}. 
This philosophical choice has a major impact on the style which is used in the remainder of the 
development. 
The undeniable existence of such an impact is perhaps one of the best arguments for why
it is interesting to try to develop mathematics from another, more constructive viewpoint. We do not
by any means want to suggest that this reasoning is faulty, on the contrary the present work could well
serve as a good example of the effect that using the axiom of choice ``early and often''
has on the subsequent style. On the other hand, the fact that there are potentially good reasons for
pursuing things from a constructive viewpoint, leaves open the potential that it might also be useful to
pursue things from a non-constructive viewpoint, and this latter is what we are doing here. 

\numero{Mathematical objects as sets with additional structure}

One of the main questions which I have tried to address is that of notation. 
In Bourbaki and thereafter, phrases of the form {\em ``un ensemble muni de \ldots ''} occur profusely.
To machine-interpret them, one needs a precise definition of the verb {\em ``munir de''}.
In this paragraph we discuss the general organizational situation this engenders. Close examination of
the question yields the conclusion that almost all of the mathematical objects we will want to consider,
can be viewed as sets of elements together with additional structure. We call the set of elements, the
{\em underlying set} of an object \verb}a}, and denote it by \verb}U a} (following the principle that notation
which occurs the most often should be as short as possible). As a first approximation, a mathematical
object could be realized as a pair \verb}a=pair x y} (in particular, \verb}a} is itself a set but we don't often
want to look at the elements of \verb}a}). Here the underlying set would be \verb}x}, and the additional structure would
be encoded in some way in the set \verb}y}. The underlying set function would just be the first projection
\verb}(U a):= (pr1 a)}.  In fact it turns out to be convenient to adopt a slightly more general strategy, 
which we discuss in the next paragraph, but for the present paragraph the reader can assume that this is
what we have done.  The main feature is that the underlying set function \verb}U} will be defined in the same
way for all different structures. 

As an example, the notion of ring would have operations \verb}plus} and \verb}times} encoded into the structure \verb}y},
so we can suppose that we have defined functions which will be written \verb}(plus a u v)} and \verb}(times a u v)}
for \verb}u,v} elements of the underlying set, i.e. \verb}inc u (U a)} and \verb}inc v (U a)}. 
Now a morphism between rings is simply a morphism of sets \verb}(U a) -> (U b)} satisfying compatibility
conditions with respect to these operations. This means that we can set up a general notation for morphisms
as being essentially triples \verb}(a,b,f)} consisting of two objects and a map \verb}f} from \verb}(U a)} to \verb}(U b)}.
In particular we can treat composition, the identity and so on, once for all possible structures.

Although we don't get this far in the present development, it is hoped that this will facilitate the
treatment of various different notions of presheaf and sheaf. As is well-known to anybody who has tried
to teach a course on sheaf theory, the notions of sheaf of sets, sheaf of groups, sheaf of rings, sheaf
of modules and so forth require at the same time a unified treatment for most things but 
specific treatments for certain things. Furthermore we would like to be able to consider automatically
a sheaf of modules as also being a sheaf of abelian groups, etc. 

To add to the complication, we would also like to look at things like simplicial sets, which are presheaves
over the category of simplices; but on the other hand we would like to look at presheaves of simplicial sets
over some other category. In particular, we need our presheaves themselves to be considered as sets 
with extra structure. This is possible to do, by taking as underlying set the disjoint union of the 
underlying sets of the values of the presheaf over the objects of the category. 

With these constructions and considerations in mind,
it seems reasonable to embark on the problem of establishing notation with the idea that we want to look
at our mathematical objects as being like pairs of the underlying set \verb}U a} plus the extra structure.

In the present development, this machinery is used principally for the case of ordered sets.
In that case, we obtain notations \verb}leq a u v} and \verb}lt a u v} for the less-than-or-equal 
(resp. strictly-less-than) relations, for \verb}inc u (U a)} and \verb}inc v (U a)}.  This works well, at least in 
the sense that it doesn't create any unexpected hitches. 

\numero{Notation}

In order to refine the approach sketched in the previous paragraph, we need to look ahead to 
the kind of constraint that we are likely to want to impose. 
The main constraint is extendability. Most basic types of objects will at some point fit into
more advanced objects with other operations. For example, a topological space is likely to be
a metric space too at some point.

When extending notation we wish to preserve the original operations, so that for example the axioms
for a metric space include as sub-properties, the axioms for a topological space (and hence we can
directly apply any theorems about topological spaces). The strategy is to try to insure a maximum of
inheritance for as low a cost as possible.

To run with the solution sketched in the previous subsection, we need to know how to encode various
operations into the structure variable which might (as a first approximation) be the second variable of a pair.
For example a ring will have structure operations \verb}plus} and \verb}times}, both of arity $2$. A topological space
will have a structure operation \verb}opens} of arity $0$. We might at some point want to mix these together.
For that reason it is inconvienent to consider the structure as being given by an ordered tuple with
different operations occupying certain places. Much more convenient is the solution which was hinted
at in Lamport-Paulson \cite{LamportPaulson}, to consider structure as being a function whose domain is
a collection of strings. That way when we add future operations there is no worry about whether 
the place they correspond to in the structure vector is already taken up by a previous operation
The places are distinguished by their strings which are possibly abbreviated versions of the name of the operation.

With this point of view, it is natural to go back on the declaration of the previous subsection, that objects
were pairs \verb}a=(x,y)} with the first projection \verb}x} being the ``underlying set''. Instead, the underlying-set
function simply becomes one of the attributes of the structure. 

We implement this. To obtain strings, we define an inductive type \verb}nota} with
$26$ constructors \verb}a_, b_, ... :nota -> nota} plus a
terminal construction \verb}DOT : nat-> nota}. The \verb}nat} variable in \verb}DOT} enables us
to specify an ``arity'' for each string (since we are interested in defining structural parameters which will
correspond to arities).

We could take a time-out here to express the need for a \verb}String} datatype in Coq with a syntactic facility to
let us write \verb}"toto"} for a string \verb}toto}---whereas in the present file an operation named ``toto'' with
arity $2$ would be written as \verb}(t_(o_(t_(o_(DOT 2)))))}. 

Thus our objects are
functions where the domain is a subset of a fixed set of ``tags'' denoting the various types of operations which
we want to have. To take an example in algebra where we would like to have rings, modules and algebras, 
there is a tag \verb}Underlying} for the underlying-set function; a tag \verb}Plus} for the
plus function; a tag \verb}Times} for the times function; and a tag \verb}Mult} for the scalar multiplication function. Respectively,
rings, modules and algebras are implemented as functions whose domains are
$$
Dom_{Ring} := \{ \verb}Underlying}, \verb}Plus}, \verb}Times} \},
$$
$$
Dom_{Module} := \{ \verb}Underlying}, \verb}Plus}, \verb}Mult} \},
$$
and
$$
Dom_{Algebra} := \{ \verb}Underlying}, \verb}Plus}, \verb}Times}, \verb}Mult} \}.
$$
In this way, the various operations are obtained by simply evaluating at the corresponding tags, so the 
``mult'' functions for modules and algebras are identical and can be assigned a single name. As a first
approximation we can write:
{\small
\begin{verbatim}
mult a x y :=  (Function.ev (Function.ev (Function.ev a Mult) x) y).
\end{verbatim}
}
In the actual file we take a somewhat more abstract approach and give a general encoding for multivariable
functions (inspired by Capretta \cite{Capretta}) which formalizes the sequence of three occurences of \verb}Function.ev} 
in the above, but
for commodity the end result is a slightly different function (this doesn't really matter though).
An explanation of that mechanism would go beyond the scope of the present note and wouldn't be very interesting
or useful, so the reader is referred directly to the source file.  Perhaps the only interesting point to note
is that we make use of the inductive structure of the datatype \verb}nat} in order to define \verb}n}-variable functions for
any \verb}n:nat}. This is an example of the benefit obtained by mixing type theory and set theory. 

Note here that the lack of type constraints and the axiom of choice come into play, because they allow us to
define function evaluation on any element as discussed previously.

A systematic application of the above principle allows a considerable amount of notational simplification.
It also allows another kind of simplification: all objects will use the same tag \verb}Underlying} for their
underlying-set functions. This allows us to develop the theory of morphisms between underlying sets, independantly
once and for all.  This is done in the module \verb}Umorphism}. The function assigning to an object \verb}a} its underlying
set is denoted \verb}a} $\mapsto$ \verb}(U a)}, with \verb}(U a)} being technically a choice of evaluation of \verb}a} 
on the element \verb}Underlying}, thus:
\begin{verbatim}
U:= fun a:E =>(Function.ev a Underlying).
\end{verbatim}
A Umorphism is a triple, or rather an object whose
domain consists of three tags \verb}Source}, \verb}Target} and \verb}Mapping}. The axioms on such a triple \verb}f} are
that \verb}(mapping f)} is a function whose domain is \verb}(U (source f))} 
and whose range is contained in \verb}(U (target f))},
where \verb}U} (resp. \verb}source}, \verb}target})
denotes the function of evaluation on the tag \verb}Underlying} (resp. \verb}Source}, \verb}Target}).  
We can define composition, identities,
inverses (when they exist), inclusions, and various lemmas about these operations, for this general notion.
Which all then applies directly to the theory of morphisms for objects whenever the morphisms can be 
faithfully expressed as maps between the underlying sets (which is really very often the case). Morphisms between
such objects are Umorphisms which are subject to additional conditions, but the constructions such as
composition and so forth are the same as in general.

We now come to another place where a mixture of type theory (with inductive definitions) and set theory
can help.  In terms of the above discussion, there is a choice to be made as to how the ``tags'' are to
be implemented as actual sets. One good choice seems to be the one 
suggested in Lamport-Paulson \cite{LamportPaulson}, that the tags be strings. Let's see why one comes to this
conclusion. 

The choice of implementation is not totally anodine, because in any given object, the tags used for
different pieces of the structure have to be distinct. The first and most obvious solution would have been to 
assign differing integer numbers to the different tags; this leads to a sort of ``Dewey decimal system'',
where we might assign for example 20-29 to category theory, 30-39 to algebra, 40-49 to topology and so forth.
This unfortunately leads to a quadratic collection of proof obligations: if we have $n$ different tags in play,
then we need the $n(n-1)/2$ statements that they are pairwise distinct.  This turns out to be a lot of different
silly things to prove.  A big improvement is obtained if we rely on \verb+Coq+'s inductive types: we can for example 
just define an inductive type containing all of the different tags as constructors.  The statement that the
constructors are different is the ``Discriminate'' tactic (see also \cite{McKinna}).  In fact we don't even need to rely on this tactic because the
difference between the constructors is integrated automatically into the inductive proof of the statements where
we need it, basically the statement that if you construct an object and then destruct it you get back what
you put in. 

The only problem with the approach outlined up until now 
was that it required putting in all of the tags which one wants to use, at the
point where the notational inductive type is defined. This results in 
an environment like ``C'' where you have to declare your data names at the beginning of the development.

An improvement, getting rid of this problem, is to define the notational inductive type (basically
a string but also with the arity included)
\verb}nota} with $27$ constructors, first one for each letter of the alphabet (26 constructors of type \verb}nota -> nota})
and a last constructor denoted \verb}dot : nat -> nota}.  Then we literally spell out the elements of \verb}nota} that
we want to use, with the \verb}dot} notation giving the arity of the operation. For example
\begin{verbatim}
Underlying := (u_(n_(d_(r_(l_(DOT 0)))))).
\end{verbatim}
This specifies \verb}Underlying} as a tag (i.e. an element of \verb}nota}) corresponding to an operation of arity $0$.
Similarly, 
\begin{verbatim}
Plus := (p_(l_(u_(s_(DOT 2))))).
\end{verbatim}
Here the arity is $2$. 
Note that it is a good idea to use rather short abbreviations because otherwise the Cases constructions 
which distinguish between different notations are expanded into huge monsters.  
With strings, anybody can add new tags (or rather, the pool of tags is free---with one generator 
for each arity---over the infinite semigroup on 
$26$ generators, so there is always room to look at new tags).

With this notational system we make a first stab, in the file \verb}notation.v}, at treating some standard
mathematical objects. In this case, ordered sets. At the end of the file \verb}ordinal.v}
we obtain a proof of Zorn's lemma. It should again
be stressed that this is not new, for example it is done in Isabelle \cite{Isabelle}, Metamath \cite{Metamath},
Mizar \cite{Mizar} and certainly other places. Rather the point is to
try out our notational style to see if it is reasonable. My conclusion was that there is probably lots of room for
improvement, but that one can at least imagine getting to more modern stuff using the system. 

While this effort certainly doesn't represent the ultimate optimal solution to the problem of fixing adequate
notation, it nonetheless should serve to indicate where the main problems are, and at least indicate that
there probably exist solutions. In this sense it should increase our confidence that the stumbling block of
notation isn't an impossible hurdle, and that it  should be possible to go on to other things.

\numero{Well-ordered sets}
One of the main new pieces of the development, with respect to where it was in \cite{ctpm},
is the file \verb}transfinite.v} in which we develop the machinery necessary to treat ordinals.
The idea is to consider a general function \verb}f:E->E} and to consider well-ordered sets \verb}a}
compatible in the sense that for any \verb}x} with \verb}inc x (U a)} we have that \verb}x=f(u)} where \verb}u} is the punctured
downward subset of \verb}x} for the order \verb}a}.  
Then we say that \verb}leq_gen f x y} if there exists a compatible \verb}a} with elements \verb}x,y} $\in $ \verb}a}
such that \verb}leq a x y}. Taking \verb}f=ID} to be the identity will yield the
order relation for ordinals (and indeed \verb}o} is an ordinal if and only if \verb}leq_gen ID o o}).

This results in essentially redoing much of the same type of argumentation which was used in the proof
of Zorn's lemma, but from a somewhat more general standpoint. 
This brings up an interesting point which illustrates something important about computer-verified
theorem proving: there is absolutely no incentive to go back and reprove Zorn's lemma using the new notation and
techniques.
Once we have given a proof, even if the proof techniques are later supplanted by better ones, 
there is no reason to redo the old proof since the only question is existence of at least one proof, which has
already been established. This is in marked contrast with the more usual situation in human-verified proofs
where it is important to have a proof which is as slick and understandable as possible, so it will be easiest to
have it verified (either by the author himself or by another reader). In that case there was an advantage in
going back and redoing the old proofs with the newer techniques. Those days are gone!

\numero{Ordinals and cardinals}

The general transfinite induction machinery of the file \verb}transfinite.v} is then used in \verb}ordinal.v} to develop
the theory of ordinals, and specially the functions \verb}wo_avatar} and \verb}wo_ordinal} which relate
an arbitrary well-ordered set \verb}a} to an ordinal. The function \verb}x}  $\mapsto$ \verb}wo_avatar a x} is the 
strictly increasing downward saturated function from the underlying set \verb}(U a)} of \verb}a}, to the 
ordinals. The ordinal \verb}(wo_ordinal a)} is the image of this function. 

The main lemma to be used later in the proof of Bernstein-Cantor-Schroeder,
is the statement \verb}suborder_wo_ordinal_decreasing}
that if \verb}sub u (U a)} then the ordinal \verb}(wo_ordinal (suborder u a))}
associated to the suborder on \verb}u}, is less than or equal to \verb}(wo_ordinal a)}.

In \verb}cardinal.v} 
we define the cardinality of a set \verb}x} to be the smallest ordinal isomorphic to \verb}x}.
This has to be proven to exist, and it needs to be shown that \verb}x} is isomorphic to \verb}(cardinality x)},
and also that isomorphic sets have the same cardinality. 
It is also convenient to have on hand the well-ordering of \verb}x} induced by a choice of such isomorphism.
Then applying the lemma described in the previous paragraph we get that if \verb}sub y x} then
\verb}ordinal_leq (cardinality y) (cardinality x)}. An immediate corollary is the Bernstein-Cantor-Schroeder theorem
which says that if \verb}x} is isomorphic to a subset of \verb}y} and \verb}y} to a subset of \verb}x}, then
\verb}(cardinality x)=(cardinality y)} and in particular \verb}x} and \verb}y} are isomorphic. It is interesting to
note that we have a totally non-constructive proof of this theorem obtained with a heavy dose of the
axiom of choice which takes care of all of the reordering of things necessary to obtain the isomorphism.
Contrast this with the classical constructive proof such as contained in \cite{HerbelinSchroeder}.
In our situation, we need the machinery of ordinal numbers and transfinite induction in any case for 
a wide range of future situations (the small object argument, for example). 
In this case the BCS theorem comes out almost for free: we avoid having to implement the constructive argument
which, although easy, nonetheless requires a certain amount of notational work. 

The philosphy at work here is that we feel that we are going to need the axiom of choice anyway
(if for no other reason than to have maximal ideals in our rings, see L. Chicli's development 
of scheme theory which starts out constructively but where he calls on existence of maximal ideals
at some point \cite{Chicli}). Consequently, why not use it as much as we like if that will 
lead to efficiencies in the proof 
development. One might well ask what is being measured when we speak of efficiency.
This is a subtle question whose response is likely to differ from one mathematician to the next.
In my own point of view, the ``cost'' is the energy it takes to plan out a statement and its proof and the
time and energy required to type it into the machine. The ``benefit'' is the contribution of that statement
to the further development (or the future development which we think we are going to do later).

We absolutely deny any originality to our approach to well-ordering, Zorn, transfinite
induction, ordinals and cardinals and the proof of BCS. I learned about all of that in the
first few weeks of Andrew Gleason's first-year graduate real analysis course.  Mathematically speaking
this is such an old subject that somebody must have exposed things in the same way, and 
as pointed out above there are many ongoing developments in other proof assistants which cover the same
material and undoubtedly overlap to whatever extent (which might go all the way to complete inclusion). 
I haven't made an extensive investigation of this issue, and it doesn't seem really necessary to do so
since the reader will easily do better on the web.

\numero{Potential additional axioms}
For the record we mention here a few possible additional axioms which might be useful. We didn't find any
particular use for them in the present development, but this might be a question of style or they might be
useful in the future.

One natural axiom would be to say that the function types \verb}x->y} or even \verb}forall a:x, y a},
should be equal to the sets of functions from \verb}x} to \verb}y} or functions on \verb}x} taking values in the
family \verb}y}. 

Another question is what to do about sub-universes of \verb}E}. The first example is \verb}Prop}.
By the axioms we have already, \verb}False} is realized as the empty set, and \verb}True} as a set with one element.
It might be good to axiomatize \verb}True = 1} so that its unique element is the emptyset. It might also be
useful to add the axiom that for \verb}X:Prop} the realization \verb}(R X)} is equal to \verb}X} (note that \verb}X}
is also an element of \verb}E} because of the universe cumulativity \verb}Prop \subset E}). We didn't encounter 
a place where this is necessary, though. If one wanted to use \verb}E':=Type :E}, say to give an inaccessible 
cardinal or for a Grothendieck universe, the corresponding axiom would probably be essential. 

\numero{Next steps}
Here is a small outline of what the next logical steps in pursuing this development might be.
The \verb}cardinal.v} file needs to be completed with basic results about finite sets, infinite sets,
also with the relationship between \verb}nat} and the finite cardinals (in the axioms we included two
axioms which imply that \verb}nat} is equal to the first infinite ordinal; thus the cardinality of a finite set
will be the element of \verb}nat} representing how many elements it has). This is not done in the current version
because of a slight lack of time. Probably the main result about finite sets is that a permutation of
a finite set can be expressed as a product of transpositions. This will be important later in the algebra
development for showing the order-invariance of a finite product of commuting elements. 

After the basic set-theory reaches a semi-satisfying first plateau, 
we should be able to access fairly rapidly the subjects of basic algebra (monoids, groups, rings, modules);
general topology, and small category theory. A more subtle question is what to do about big categories.
The Coq proof environment technically could provide the tools necessary to introduce the notion of
class (by going to a variable \verb}Type_i} higher than the one used for \verb}E}), but this seems like it
would be complicated, and would put us in the position of wanting to have a refined ``typical ambiguity''
like what Feferman suggests in \cite{Feferman} which is a difficult research subject.
On the other hand, our notion of ``underlying-set function'' implemented with the function \verb}U} and
the module \verb}Umorphism} opens up the possibility of considering a limited kind of
big category which would be a category whose objects are ``sets with extra structure'',
that is sets \verb}x} thinking of the underlying set being \verb}(U x)}; the morphisms of the category would
be certain \verb}Umorphisms}. Since pretty much all of the objects we are interested in can be 
worked into this mold, and since it is probably fairly easy to rapidly develop the small part of
general machinery common to all the examples, this is probably an efficient solution.
This choice would guide our implementation of the notion of presheaf: a presheaf should be thought of
as a functor from a small category to one of these ``Ucategories''. This can be made explicit,
and furthermore such a presheaf can itself be considered as an underlying set with extra structure,
where the underlying set is the set of ordered pairs \verb}(o,x)} where \verb}o} is an object of the base category
and \verb}x} is an element of the underlying set of the value of the presheaf on \verb}o}. Thus given a \verb}Ucategory} \verb}F} 
and a small category \verb}c} we can define a new \verb}Ucategory} called \verb}Presh c F}.  This would formalize
the presheaf construction in a way adapted to things like simplicial homotopy theory
and geometry of ringed spaces. 

The possibility of continuing the present development along the above lines, would be a crucial test of
its potential validity as a groundwork for reaching modern research-level mathematics 
with a proof verification tool.

 \end{document}